\documentclass[12pt,reqno]{amsart}
\usepackage{amsmath, amsthm, amssymb, stmaryrd}
\usepackage{hyperref}
\usepackage{enumerate}
\usepackage{url}

\let\OLDthebibliography\thebibliography
\renewcommand\thebibliography[1]{
  \OLDthebibliography{#1}
  \setlength{\parskip}{0pt}
  \setlength{\itemsep}{0pt plus 0.3ex}
}

\usepackage{mathtools}

\usepackage{multirow, array}
\usepackage{placeins}

\usepackage{caption} 
\advance \topmargin by -\headheight
\advance \topmargin by -\headsep
     
\evensidemargin \oddsidemargin
\newtheorem{thm}{Theorem}[section]

\theoremstyle{definition}

\theoremstyle{remark}

\numberwithin{equation}{section}

\newcommand*\wrapletters[1]{\wr@pletters#1\@nil}
\def\wr@pletters#1#2\@nil{#1\allowbreak\if&#2&\else\wr@pletters#2\@nil\fi}

\def\alp{{\alpha}}

\def\del{{\delta}}

\def\eps{\varepsilon}

\def\le{\leqslant} \def\ge{\geqslant}

\def \bN {\mathbb N}

\def \bR {\mathbb R}
\def \bZ {\mathbb Z}

\def \balp {\boldsymbol{\alp}}

\def \cC {\mathcal C}

\begin{document}
\title[A note on rational points near planar curves]{A note on rational points near planar curves}
\author[Sam Chow]{Sam Chow}
\address{Department of Mathematics, University of York,
Heslington, York, YO10 5DD, United Kingdom}
\email{sam.chow@york.ac.uk}
\thanks{}
\date{}
\begin{abstract} Under fairly natural assumptions, Huang counted the number of rational points lying close to an arc of a planar curve. He obtained upper and lower bounds of the correct order of magnitude, and conjectured an asymptotic formula. In this note, we establish the conjectured asymptotic formula.
\end{abstract}
\maketitle

\section{Introduction}
\label{intro}

Let $f$ be a real-valued function defined on a compact interval $I = [\rho, \xi] \subseteq \bR$. For positive real numbers $\del \le 1/2$ and $Q \ge 1$, define
\[
\tilde{N}_f(Q, \del) = \# \left \{ 
  \begin{aligned}
 (a,b,q) \in \bZ^3: \quad &1 \le q \le Q, a/q \in I, \gcd(a,b,q)=1,\\
  & |f(a/q)-b/q| < \del/Q 
  \end{aligned}
\right\}.
\]
Roughly speaking, this counts the number of rational points with denominator at most $Q$ that lie within $\del Q^{-1}$ of the curve $\cC_f = \{ (x, f(x)): x \in I \}$. Huang \cite[Theorem 2]{Hua2015} estimated this quantity. As discussed in \cite{Hua2015}, such estimates are readily applied to the Lebesgue theory of metric diophantine approximation.

\begin{thm}[Huang] \label{HuangThm} Let $0 < c_1 \le c_2$. Assume that $f: I \to \bR$ is a $C^2$ function satisfying 
\[
c_1 \le |f''(x)| \le c_2 \qquad (x \in I),
\]
with Lipschitz second derivative. Assume further that 
\begin{equation} \label{hyp}
1/2 \ge \del > Q^{\eps - 1},
\end{equation}
for some $\eps \in (0,1)$. Then
\begin{equation} \label{Huang2}
\frac{2 \sqrt3}{9 \zeta(3)} + O(Q^{-\eps/2}) \le \frac{\tilde{N}_f(Q, \del)}{|I| \del Q^2} \le \frac1{\zeta(3)} + O(Q^{-\eps/2}).
\end{equation}
The implied constant depends on $I, c_1, c_2, \eps$ and the Lipschitz constant; it is independent of $f, \del$ and $Q$. 
\end{thm}

Theorem \ref{HuangThm} sharpened the upper bounds obtained by Huxley \cite{Hux1994} and Vaughan--Velani \cite{VV2006}, as well as the lower bounds obtained by Beresnevich--Dickinson--Velani \cite{BDV2007} and Beresnevich--Zorin \cite{BZ2010}.

The purpose of this note is to squeeze together the constants in \eqref{Huang2}, so as to confirm Huang's conjectured asymptotic formula
\begin{equation} \label{asymp}
\tilde{N}_f(Q,\del) \sim \frac2{3 \zeta(3)} |I| \del Q^2 \qquad (Q \to \infty),
\end{equation}
within the range \eqref{hyp}. The asymptotic formula \eqref{asymp} follows straightforwardly from our theorem, which we state below and establish in the next section.

\begin{thm} \label{MainTheorem} Assume the hypotheses of Theorem \ref{HuangThm}. Let $\eta > 0$ and
\[
0 < \tau < \eps/2.
\]
Then
\[
\frac2{3 \zeta(3)} - \eta + O(Q^{-\tau}) \le \frac{\tilde{N}_f(Q,\del)}{|I| \del Q^2} \le \frac2{3 \zeta(3)} + \eta + O(Q^{-\tau}).
\]
The implied constant depends on $I, c_1, c_2, \eps, \eta$ and the Lipschitz constant.
\end{thm}

We use Landau and Vinogradov notation: for functions $f$ and positive-valued functions $g$, we write $f \ll g$ or $f = O(g)$ if there exists a constant $C$ such that $|f(x)| \le C g(x)$ for all $x$. If $S$ is a set, we denote the cardinality of $S$ by $\# S$. 

The author is supported by EPSRC Programme Grant EP/J018260/1, and thanks Faustin Adiceam for a discussion.

\section{The count}
\label{proof}

In this section, we prove Theorem \ref{MainTheorem}. For positive real numbers $\del \le 1/2$ and $Q \ge 1$, define the auxiliary counting function
\[
\hat{N}_f(Q, \del) = \# \left \{ 
  \begin{aligned}
 (a,b,q) \in \bZ^3: \quad &1 \le q \le Q, a/q \in I, 
 \\ & 
  \gcd(a,b,q)=1,|f(a/q)-b/q| < \del/q
  \end{aligned}
\right\}.
\]
With the same assumptions as in Theorem \ref{HuangThm}, Huang \cite[Corollary 1]{Hua2015} showed that
\begin{equation} \label{Huang1}
\hat{N}_f(Q, \del) = (\zeta(3)^{-1} + O(Q^{-\eps/2})) \cdot |I| \del Q^2.
\end{equation}
Let $t \in \bN$, $1/2 < \alp < 1$ and
\[
\alp_i = \alp^i \qquad (0 \le i \le t).
\]
We will have $t \ll_\eta 1$, so the hypothesis \eqref{hyp} is satisfied with $2 \tau$ in place of $\eps$ and $(\alp_iQ, \alp_j \del)$ in place of $(Q,\del)$, whenever $Q$ is large and $0 \le i,j \le t$. In particular \eqref{Huang1} holds with these adjustments, so
\begin{equation} \label{Huang}
\hat{N}_f(\alp_i Q, \alp_j \del) = \Bigl(\frac{\alp_i^2 \alp_j}{\zeta(3)} + O(Q^{-\tau}) \Bigr) \cdot |I| \del Q^2 \qquad (0 \le i,j \le t).
\end{equation}

Employing \eqref{Huang}, we have
\begin{align*}
\tilde{N}_f(Q, \del) &\ge
\sum_{i=1}^t \# \left \{ 
  \begin{aligned}
 (a,b,q) \in \bZ^3: \quad &\alp_i Q < q \le \alp_{i-1} Q, a/q \in I, 
 \\ &
 \gcd(a,b,q)=1, |f(a/q)-b/q| < \alp_i \del/q
  \end{aligned}
\right\} \\
&= \sum_{i=1}^t (\hat{N}_f(\alp_{i-1}Q,\alp_i\del) - \hat{N}_f(\alp_i Q,\alp_i\del)) \\
&= \sum_{i=1}^t \Bigl(\frac{\alp_{i-1}^2 \alp_i - \alp_i^3}{\zeta(3)} + O(Q^{-\tau}) \Bigr) \cdot |I| \del Q^2.
\end{align*}
Now
\begin{equation} \label{lower}
\tilde{N}_f(Q,\del) \ge \Bigl(\frac{X(\balp)}{\zeta(3)} + O(tQ^{-\tau}) \Bigr) \cdot |I| \del Q^2,
\end{equation}
where
\[
X(\balp) = \sum_{i \le t} (\alp_{i-1}^2 \alp_i - \alp_i^3).
\]
We compute that
\begin{align*}
X(\balp) &= (\alp - \alp^3) \sum_{j=0}^{t-1} (\alp^3)^j = \frac{(\alp - \alp^3) (1- \alp^{3t})}{1- \alp^3} \\
&= (1-\alp^{3t})(1- (1+\alp + \alp^2)^{-1}).
\end{align*}
Choosing $\alp$ close to 1, and then choosing $t \ll_\eta 1$ large, gives
\[
X(\balp) \ge 2/3 - \zeta(3) \eta.
\]
Substituting this into \eqref{lower} yields the desired lower bound.

We attack the upper bound in a similar fashion, but there is an extra term to consider. By \eqref{Huang}, we have
\begin{align*}
& \tilde{N}_f(Q, \del) - \tilde{N}_f(\alp_t Q, \alp_t \del) 
\\ &\le 
\sum_{i=1}^t \# \left \{ 
  \begin{aligned}
 (a,b,q) \in \bZ^3: \quad &\alp_i Q < q \le \alp_{i-1} Q, a/q \in I, \gcd(a,b,q)=1,\\
  & |f(a/q)-b/q| < \alp_{i-1} \del/q
  \end{aligned}
\right\} \\
&= \sum_{i=1}^t (\hat{N}_f(\alp_{i-1}Q,\alp_{i-1}\del) - \hat{N}_f(\alp_i Q,\alp_{i-1}\del)) \\
&= \sum_{i=1}^t \Bigl(\frac{\alp_{i-1}^3 - \alp_{i-1}\alp_i^2}{\zeta(3)} + O(Q^{-\tau}) \Bigr) \cdot |I| \del Q^2.
\end{align*}
Now 
\[
\tilde{N}_f(Q, \del) - \tilde{N}_f(\alp_t Q, \alp_t \del)  \le \Bigl(\frac{Y(\balp)}{\zeta(3)} + O(t Q^{-\tau}) \Bigr) \cdot |I| \del Q^2,
\]
where
\[
Y(\balp) = \sum_{i \le t} (\alp_{i-1}^3 - \alp_{i-1} \alp_i^2).
\]
Here
\[
Y(\balp)= \alp^{-1} X(\balp) \le \frac{1- \alp^2}{1-\alp^3} = \frac{1+\alp}{1+\alp+\alp^2}.
\]
Choosing $\alp$ close to 1 gives $Y(\balp) \le 2/3 + \zeta(3)\eta/2$, and so
\begin{equation} \label{upper}
\tilde{N}_f(Q, \del) \le \tilde{N}_f(\alp_t Q, \alp_t \del) + \Bigl(\frac2{3\zeta(3)} + \frac \eta2 + O(t Q^{-\tau}) \Bigr) \cdot |I| \del Q^2.
\end{equation}

For the first term on the right hand side of \eqref{upper}, we bootstrap Huang's upper bound \eqref{Huang2}. This gives
\[
\tilde{N}_f(\alp_t Q, \alp_t \del) \le \Bigl(\frac{\alp_t^3}{\zeta(3)} + O(Q^{-\tau})\Bigr) \cdot |I| \del Q^2.
\]
Choosing $t \ll_\eta 1$ large, so that $\alp_t^3 \le \zeta(3) \eta/2$, we now have
\[
\tilde{N}_f(\alp_t Q, \alp_t \del) \le \Bigl(\frac \eta 2 + O(Q^{-\tau}) \Bigr) \cdot |I| \del Q^2.
\]
Substituting this into \eqref{upper} provides the sought upper bound, completing the proof of the theorem.

\bibliographystyle{amsbracket}
\providecommand{\bysame}{\leavevmode\hbox to3em{\hrulefill}\thinspace}

\end{document}